\documentclass[twoside,12pt]{article}
 \usepackage{bbm}
 \usepackage{mathrsfs}
  \usepackage{amsfonts}
\usepackage{amsmath}
\newcommand{\bref}[5]{\bibitem{#1} {#2} {\it #3} {\bf #4}#5.}

\def\coker{\operatorname{coker}}

\begin{document}
  \title{The critical group of $K_m\times C_n$
 \thanks{Supported by NSF of the People's Republic of China(Grant
 No. 10671191,  No. 10301031 and No. 10871189).}
  }
 \author{Jian Wang,  Yong-Liang Pan\thanks{Corresponding author. Email: ylpan@ustc.edu.cn},\,\,
 Jun-Ming Xu
 \\
  {\small Department of Mathematics, University of Science and Technology
         of China}\\
  {\small Hefei, Auhui 230026, The People's Republic of China}\\
}
\date{}
\maketitle {\centerline{\bf\sc Abstract}\vskip 8pt In this paper,
the structure of the critical group of the graph $K_m\times C_n$ is
determined, where $m,\, n\ge 3$.



\section{Introduction and statement of results}

\indent The critical group of a connected graph is a finite abelian
group whose structure is a subtle isomorphism invariant of the
graph.
It is closely connected with the graph Laplacian.\\
\indent  Let $G=(V,E)$ be a finite connected graph without
self-loops, but with multiple edges allowed. Then the Laplacian
matrix of $G$ is the $|V|\times |V|$ matrix defined by
 $$L(G)_{uv}=\left\{
\begin{array}{ll}
d(u), &  \text{if}\hspace{0.2cm}u=v, \\
-a_{uv}, & \text{if}\hspace{0.2cm}u\not=v,
\end{array}
\right.\eqno(1.1) $$ where $a_{uv}$ is the number of the edges
joining $u$ and $v$, and $d(u)$ is the degree of $u$.\\
\indent Regarding $L(G)$ as a homomorphism ${\mathbb{Z}}^{|V|}
\rightarrow {\mathbb{Z}}^{|V|}$, its cokernel
coker$(L(G))={\mathbb{Z}}^{|V|}/\text{im}\, (L(G))$ is an abelian
group. For $1\leq i\leq |V|$, let $e_i=(0,\cdots,0,1,0,\cdots,0)^t
\in {\mathbb{Z}}^{|V|}$, be the $i$-th standard basis, and $x_i$ be its image in
$\coker(L(G))$. We know that $\coker(L(G))$ is determined by the
generators $x_1,\cdots, x_{|V|}$ and the relations  $(x_1,\cdots,x_{|V|})L(G)=0$.
Since $L(G)$ is symmetric, we can rewrite the
relations  as follows
$$\left\{\begin{array}{ll}
l_{11}x_1+l_{12}x_2+\cdots+l_{1|V|}x_{|V|}=0,\\
l_{21}x_1+l_{22}x_2+\cdots+l_{2|V|}x_{|V|}=0,\\
\vdots\\
l_{|V|1}x_1+l_{|V|2}x_2+\cdots+l_{|V||V|}x_{|V|}=0.\\
\end{array}\right.\eqno(1.2)$$
\indent Two integral matrices $A$ and $B$  are
equivalent (written $A \sim B$) if there are unimodular  matrices
$P$ and $Q$ such that $B=PAQ$ (An integral matrix $P$ is unimodular
if $P^{-1}$ is also integral, i.e., if $\det P=\pm 1$.).
Equivalently, $B$ is obtainable from $A$ by a sequence of elementary
row and column operations: (1) the interchange of two rows or
columns, (2) the multiplication of  any row or column by $-1$,
(3) the addition of any integer times of one row (resp. column) to another row (resp. column).\\
\indent It is easy to see that $A\sim B$  implies that
coker$(A)\cong$ coker$(B)$. The Smith normal form is a diagonal
canonical form for our equivalence relation: every $n\times n$ integral matrix
$A$  is equivalent to a unique diagonal matrix diag$(s_1(A),\cdots,s_n(A))$, where $s_i(A)$ divides $s_{i+1}(A)$ for
$i=1,2,\cdots, n-1$. The $i-$th diagonal entry of the Smith normal form of $A$ is usually called the $i-$th invariant factor of $A$.  We will use the fact that the values $s_i(A)$ can
also be interpreted as follows: for each $i$, the product
$s_1(A)s_2(A)\cdots s_i(A)$ is the greatest common divisor of all $i\times i$
minors  of $A$.\\
\indent  The classification theorem for finitely generated abelian
groups asserts that  coker($L(G))$ has a direct sum decomposition
$$\mbox{coker}(L(G))\cong\left({\mathbb{Z}}/ t_1 {\mathbb{\mathbb{Z}}}\right)\oplus \left({\mathbb{Z}}/t_2 {\mathbb{Z}}\right)
\oplus\cdots\oplus\left({\mathbb{Z}}/ t_{|V|}
\mathbb{Z}\right),\eqno(1.3)$$ where the nonnegative integers $t_i$  are
the diagonal entries of the Smith normal form of the relation matrix
$L(G)$, of course, they satisfy that $t_i$ divides $t_{i+1},\, (1\le i<|V|)$.
Since $G$ is connected, it is not hard to see that $L(G)$ has rank
$|V|-1$, and the kernel of $L(G)$ is spanned by the vectors in
${\mathbb{R}}^{|V|}$ which are constant on the vertices. It follows that
$t_{|V|}=0$ and
$t_1\cdots t_{|V|-1}\not=0$. \\
\indent Now we can write
$$\mbox{coker}(L(G))={\mathbb{Z}}^{|V|}/\text{im}\, (L(G))\cong {\mathbb{Z}}\oplus K(G),\eqno(1.4)$$
where
$$K(G)=\left({\mathbb{Z}}/ t_1 {\mathbb{Z}}\right)\oplus \left({\mathbb{Z}}/t_2 {\mathbb{Z}}\right)
\oplus\cdots\oplus\left({\mathbb{Z}}/ t_{|V|-1}
{\mathbb{Z}}\right).\eqno(1.5)$$
\indent The finite abelian group
$K(G)$  is defined to be the critical group of $G$. And we will call
the positive integers $t_1,\cdots, t_{|V|-1}$ the invariant factors
of $K(G)$. The critical group $K(G)$ is also known as the Picard
group and the Jacobian group of $G$ in [1,\,2,\,3 ], while in the
physics literature it is known as the abelian sandpile group, and it
has a close connection with the critical configuration in a certain
dollar game on $G$, see [3, 9]. For the general theory of the
critical group, we refer the reader to Biggs [2, 3],
 Godsil [9, Chapter 14], Cori, et al. [5, 6], Dartois et al. [8], and  Bacher, et al. [1].\\
\indent The well known Kirchhoff's Matrix-Tree Theorem [9, Theorem
13.2.1] shows that $t_1\cdots t_{|V|-1}$
 equals the number $\kappa$ of spanning trees of $G$. It follows that  the invariant factors of $K(G)$ can be used to
distinguish pairs of non-isomorphic graphs which have the same
$\kappa$, and so there is considerable interest in their properties.
If $G$ is a simple connected graph, then its Laplacian matrix $L(G)$
has some entry which is equal to $-1$. Since the invariant factor
$t_1$ of $K(G)$ is equal to the greatest common divisor of all the
entries of $L(G)$, it follows that $t_1$ must be equal to 1.
But the other invariant factors of $K(G)$  are not easy to be determined.\\
\indent Compared to the number of the results on the spanning tree
number $\kappa$, there are relatively few results describing the
critical group structure of $K(G)$ in terms of the structure of $G$.
Recently, there are some families of graphs for which the critical
group structure has been completely determined:  wheel graphs [3];
cycles [14]; complete graphs [12];
complete multipartite graphs and cartesian products of complete
graphs [11]; a subclass of the threshold graphs [4]; the M\"{o}bius
ladder graphs [7]; the Cayley graph $\mathcal{D}_n$ of the dihedral
group [8];
the square cycle graphs $C_n^2$ [10]; etc. \\
\indent Given two disjoint graphs $G_1=(V_1,E_1)$ and
$G_2=(V_2,E_2)$, the Cartesian product  of them is denoted by
$G_1\times G_2$. It has vertex set $V_1\times V_2=\{(u_i,v_j)|
u_i\in V_1, v_j\in V_2\}$, where $(u_1,v_1)$ is adjacent to
$(u_2,v_2)$ if and only if $u_1=u_2$ and $(v_1,v_2)\in E_2$, or
$(u_1,u_2)\in E_1$ and $v_1=v_2$. One may view $G_1\times G_2$ as
the graph obtained from $G_2$ by replacing each of its  vertices
with a copy of $G_1$, and each of its edges with $|V_1|$ edges
joining corresponding vertices
of $G_1$ in the two copies.\\
\setlength{\unitlength}{1cm}
\begin{picture}(18,6)(-1,0)
\put(1,0){\line(0,1){3}} \put(1,0){\line(1,1){1}}
\put(2,4){\line(-1,-1){1}} \put(2,4){\line(0,-1){3}}
\put(3,0){\line(0,1){3}} \put(3,0){\line(1,1){1}}
\put(4,4){\line(-1,-1){1}} \put(4,4){\line(0,-1){3}}
\put(5,0){\line(0,1){3}} \put(5,0){\line(1,1){1}}
\put(6,4){\line(-1,-1){1}} \put(6,4){\line(0,-1){3}}
\put(9,0){\line(0,1){3}} \put(9,0){\line(1,1){1}}
\put(10,4){\line(-1,-1){1}} \put(10,4){\line(0,-1){3}}
\put(11,0){\line(0,1){3}} \put(11,0){\line(1,1){1}}
\put(12,4){\line(-1,-1){1}} \put(12,4){\line(0,-1){3}}
\put(1,0){\line(1,0){5}}\put(6,0){\line(1,0){0.2}}\put(6.4,0){\line(1,0){0.2}}\put(6.8,0){\line(1,0){0.2}}
\put(7.2,0){\line(1,0){0.2}}\put(7.6,0){\line(1,0){0.2}}\put(8,0){\line(1,0){3}}
\put(2,1){\line(1,0){5}}\put(7,1){\line(1,0){0.2}}\put(7.4,1){\line(1,0){0.2}}\put(7.8,1){\line(1,0){0.2}}
\put(8.2,1){\line(1,0){0.2}}\put(8.6,1){\line(1,0){0.2}}\put(9,1){\line(1,0){3}}
\put(2,4){\line(1,0){5}}\put(7,4){\line(1,0){0.2}}\put(7.4,4){\line(1,0){0.2}}\put(7.8,4){\line(1,0){0.2}}
\put(8.2,4){\line(1,0){0.2}}\put(8.6,4){\line(1,0){0.2}}\put(9,4){\line(1,0){3}}
\put(1,3){\line(1,0){5}}\put(6,3){\line(1,0){0.2}}\put(6.4,3){\line(1,0){0.2}}\put(6.8,3){\line(1,0){0.2}}
\put(7.2,3){\line(1,0){0.2}}\put(7.6,3){\line(1,0){0.2}}\put(8,3){\line(1,0){3}}
\put(1,3){\line(1,-2){1}} \put(1,0){\line(1,4){1}}
\put(3,3){\line(1,-2){1}} \put(3,0){\line(1,4){1}}
\put(5,3){\line(1,-2){1}} \put(5,0){\line(1,4){1}}
\put(9,3){\line(1,-2){1}} \put(9,0){\line(1,4){1}}
\put(11,3){\line(1,-2){1}} \put(11,0){\line(1,4){1}}
\qbezier(1,0)(6,-3)(11,0) \qbezier(1,3)(6,6)(11,3)
\qbezier(2,1)(7,-2)(12,1) \qbezier(2,4)(7,7)(12,4)
\put(1,0){\oval(0.15,0.1)} \put(1,3){\oval(0.15,0.1)}
\put(2,1){\oval(0.15,0.1)} \put(2,4){\oval(0.15,0.1)}
\put(3,0){\oval(0.15,0.1)} \put(3,3){\oval(0.15,0.1)}
\put(4,1){\oval(0.15,0.1)} \put(4,4){\oval(0.15,0.1)}
\put(5,0){\oval(0.15,0.1)} \put(5,3){\oval(0.15,0.1)}
\put(6,1){\oval(0.15,0.1)} \put(6,4){\oval(0.15,0.1)}
\put(9,0){\oval(0.15,0.1)} \put(9,3){\oval(0.15,0.1)}
\put(10,1){\oval(0.15,0.1)} \put(10,4){\oval(0.15,0.1)}
\put(11,0){\oval(0.15,0.1)} \put(11,3){\oval(0.15,0.1)}
\put(12,1){\oval(0.15,0.1)} \put(12,4){\oval(0.15,0.1)}
\put(1.1,-0.3){$v_{0,0}$} \put(2,0.7){$v_{0,1}$}
\put(1.1,3.3){$v_{0,m-1}$} \put(3,-0.3){$v_{1,0}$} \put(4,0.7){$v_{1,1}$}
\put(3.1,3.3){$v_{1,m-1}$} \put(5,-0.3){$v_{2,0}$} \put(6,0.7){$v_{2,1}$}
\put(5.1,3.3){$v_{2,m-1}$} \put(9,-0.3){$v_{n-2,0}$}
\put(10,0.7){$v_{n-2,1}$}
 \put(9,3.3){$v_{n-2,m-1}$}
\put(11,-0.3){$v_{n-1,0}$} \put(12,0.7){$v_{n-1,1}$}
 \put(11,3.3){$v_{n-1,m-1}$}
\put(3.8,-2.5){Fig. 1. Graph $K_m\times C_n$.}
\end{picture}
\\
\vskip 2.5cm
\indent The structure of the critical group of $K_m\times P_n$ has
been obtained in [13], where $K_m$ is the complete graph on $m$
vertices and $P_n$ is the path on $n$ vertices. In this paper  we
will describe the structure of the critical group on $K_m\times C_n$
with $n, m\ge 3$, where  $C_n$ is the cycle on $n$ vertices. From
the definition of the Cartesian product of two graphs, it is easy to
see that there are $n$ layers of $K_m\times C_n$, each of which is a
copy of $K_m$. Let ${\mathbb{Z}}_n$ denote ${\mathbb{Z}}/n{\mathbb{Z}}$, then for $i\in
{{\mathbb{Z}}}_n,\;  j\in {\mathbb{Z}}_m$, we may let $v_{i,j}$ denote the $j$-th vertex
in the $i$-th layer of $K_m\times C_n$.  The vertex $v_{i,j}$ is
adjacent to vertices $v_{l,j}$ with
$l=i\pm 1(\mod n)$, and to the vertices $v_{i,k}$,  $k\in {\mathbb{Z}}_m, k\not=j$ (See Fig. 1).\\

\indent Before the main result can be stated, we need some technical definitions.\\
\indent If $m$ is a positive integer, let
$\alpha=\frac{1}{2}(m+2+\sqrt{m^2+4m})$,\,
$\beta=\frac{1}{2}(m+2-\sqrt{m^2+4m})$. Then for $p\in \mathbb{Z}$, we
set
$u_p:=\frac{1}{\alpha-\beta}\left(\alpha^p-\beta^p\right),$\,$v_p:=\alpha^p+\beta^p,$\,
$\tau_p:=\frac{1}{m}(p-u_p)$, $h_p:=u_p+u_{p+1}$, and
$g_p:=\tau_p+\tau_{p+1}$. For the integers
$a_1,\,a_2,\cdots,a_k$,  we will let $(a_1,a_2,\cdots, a_k)$ denote their greatest common divisor,
and  use $a_1\mid a_2\mid\cdots\mid a_k$  to mean that $a_1$ divides $a_2$, $a_2$ divides $a_3$, etc.\\

\indent Now, we can state our main result in this article as follows.\\
\noindent{\bf Theorem 1.1}\,  If $n=2s+1$, the critical group of
$K_m\times C_n$ $(m,\, n\ge 3)$ is
$$\mathbb{Z}_{(n,g_s)}\oplus \mathbb{Z}_{h_s}\oplus \underset{m-2}{\underbrace{\mathbb{Z}_{h_s}\oplus
\cdots\oplus \mathbb{Z}_{h_s}}}\oplus \mathbb{Z}_{\gamma}\oplus
\underset{m-3}{\underbrace{\mathbb{Z}_{mh_s}\oplus\cdots\oplus
\mathbb{Z}_{mh_s}}}\oplus \mathbb{Z}_{\varphi},$$ where $\gamma=\frac{h_s}{(n,g_s)}(n, h_s)$ and
$\varphi=\frac{nmh_s}{(n,h_s)}$.\\
\indent If $n=2s$, the critical group of $K_m\times C_n$ $(m,\, n\ge
3)$ is
$$\mathbb{Z}_{(u_s, 2\tau_s)}\oplus \mathbb{Z}_{\zeta} \oplus \underset{m-3}{\underbrace{\mathbb{Z}_{(m, 2)u_s}
\oplus\cdots\oplus \mathbb{Z}_{(m, 2)u_s}}}\oplus \mathbb{Z}_{\eta}\oplus \mathbb{Z}_{\rho}\oplus
\underset{m-3}{\underbrace{\mathbb{Z}_{\chi}\oplus\cdots\oplus
\mathbb{Z}_{\chi}}}\oplus \mathbb{Z}_{\xi}, $$ where $$\left\{\begin{array}{cl}
\zeta &=\frac{u_s\left(n,\, u_s,\, 4\tau_s\right)}{(u_s, 2\tau_s)},\\
\eta  &=\frac{u_s(m,2) \left(n,  u_s-4\tau_s\right)}
{\left(n,\, u_s,\, 4\tau_s\right)},\\
\rho  &=\frac{(m+4)u_s\left( mn,\, (m+4)u_s,\, 2n\right)}{\left(n,\, u_s-4\tau_s\right)(m,\, 2)},\\
\chi&=\frac{m(m+4)u_s}{\left(m, 2\right)},\\
\xi   &=\frac{nm(m+4)u_s}{\left(mn,\ (m+4)u_s,\ 2n\right)}.
\end{array}\right.$$\\

\indent An immediate consequence of theorem 1.1 is the following Corollary.\\
\noindent{\bf Corollary 1.2}\, The spanning tree number of
$K_m\times C_n$ is
$$\frac{n}{m}\left(\left(\frac{m+2+\sqrt{m^2+4m}}{2}\right)^n+\left(\frac{m+2-\sqrt{m^2+4m}}{2}\right)^n-2\right)^{m-1}.$$

\section{Propositions and Lemmas}

\indent We  first present some obvious and some less
obvious Propositions of the sequences
$u_p$, $v_p$, $\tau_p$, $h_p$ and $g_p$.\\
\indent  Note that
$\alpha^p\mp\beta^p=(\alpha+\beta)(\alpha^{p-1}\mp\beta^{p-1})-\alpha\beta(\alpha^{p-2}\mp\beta^{p-2})$.
With the above definitions, it is easy to see that
$\alpha+\beta=m+2$ and $\alpha\beta=1$.
So we have the following Proposition 2.1.\\
\noindent{\bf Proposition 2.1} If $p$ is integral, then
$$\left\{\begin{array}{ll}
u_p=(m+2)u_{p-1}-u_{p-2},\\
u_0=0,\quad \quad u_1=1,\\
v_p=(m+2)v_{p-1}-v_{p-2},\\
v_0=2,\quad v_1=m+2.
\end{array}\right.\eqno(2.1)$$
\indent From $(2.1)$, it is easy to see that for every integer $p$,
$u_p$ and $v_p$ are integral.
The following Propositions 2.2 can be proved by induction on $p$.\\
\noindent{\bf Proposition 2.2} If $p$ is integral, then
$$u_p\equiv p\ (\mod\ m),\  \ v_p\equiv 2\ (\mod\ m). \eqno(2.2)$$
\indent By (2.2), we see that $m\mid (p-u_p)$, i.e., $\tau_p$ is
integral for $p\in \mathbb{Z}$.
In fact, we further have the following Proposition (2.3).\\
\noindent{\bf Proposition 2.3} If $p$ is integral, then
$$\tau_p=(m+2)\tau_{p-1}-\tau_{p-2}-(p-1).\eqno(2.3)$$
\noindent{\bf Proof}\, Since $u_p=p-m\tau_p$, it follows from Proposition 2.1 that
$p-m\tau_p=(m+2)(p-1-m\tau_{p-1})-(p-2-m\tau_{p-2})$.
So $m\tau_p=m(m+2)\tau_{p-1}-m\tau_{p-2}-m(p-1)$ and then (2.3) holds.\hfill$\Box$\\
\noindent{\bf Proposition 2.4} If $p$ is a nonnegative integer, then
$$\begin{array}{ll}
&u_{p-1}u_{p+1}-u_p^2-1+(u_{p+1}-u_{p-1})\\
&=v_p-2=\left\{
\begin{array}{ll}
mh_s^2, &\mbox{if}\; p=2s+1,\\
m(m+4)u_s^2, &\mbox{if}\; p=2s.
\end{array}
\right. \end{array}\eqno (2.4)$$
\noindent{\bf Proof}\, A direct calculation can show $$
\begin{array}{ll}&u_{p-1}u_{p+1}-u_p^2-1+(u_{p+1}-u_{p-1})\\
&=\frac{\alpha^{p-1}-\beta^{p-1}}{\alpha-\beta}\cdot\frac{\alpha^{p+1}-\beta^{p+1}}{\alpha-\beta}
-\left(\frac{\alpha^{p}-\beta^{p}}{\alpha-\beta}\right)^2-1
+\left(\frac{\alpha^{p+1}-\beta^{p+1}}{\alpha-\beta}-\frac{\alpha^{p-1}-\beta^{p-1}}{\alpha-\beta}\right)\\
&=-1-1+\alpha^p+\beta^p=v_p-2.\end{array}$$
So the first equality of (2.4) holds. Now we verify the second equality.\\
\indent If $p=2s+1$, then
$$\begin{array}{ll}mh_s^2&=m(u_{s+1}+u_s)^2=\frac{m}{(\alpha-\beta)^2}(\alpha^{s+1}-\beta^{s+1}+\alpha^s-\beta^s)^2\\
&=\frac{m}{m^2+4m}(\alpha^{2s+2}+\beta^{2s+2}+\alpha^{2s}+\beta^{2s}+2\alpha^{2s+1}+2\beta^{2s+1}-2-2-2\alpha-2\beta)\\
&=\frac{1}{m+4}((m+4)v_p-2(m+4))=v_p-2.\end{array}$$
\indent If $p=2s$, then
$$m(m+4)u_s^2=\frac{m(m+4)}{(\alpha-\beta)^2}(\alpha^s-\beta^s)^2
           =\alpha^{2s}+\beta^{2s}-2=v_p-2.$$
\hfill$\Box$\\
\noindent{\bf Proposition 2.5} If $p$ is integeral, then
$$(u_{p+1}-1,\ u_p)=(u_p,\ u_{p-1}+1)=
\left\{\begin{array}{ll}
h_s, &\mbox{if}\; p=2s+1,\\
(m,\ 2)u_s, &\mbox{if}\; p=2s.
\end{array}
\right. \eqno(2.5)
$$
\noindent{\bf Proof} For $i\in\mathbb{Z}$, set $\theta_i:=u_{p-i}+u_i$.
Note that $\alpha\beta=1$ implies that $u_{-i}=-u_i$, then it
follows from (2.1) that
$$\begin{array}{ll}
\theta_{i+1}&=u_{p-i-1}+u_{i+1}=-u_{i+1-p}+u_{i+1}\\
&=-((m+2)u_{i-p}-u_{i-1-p})+(m+2)u_i-u_{i-1}\\
&=(m+2)u_{p-i}+u_{i-1-p}+(m+2)u_i-u_{i-1}\\
&=(m+2)(u_{p-i}+u_i)-(u_{p-(i-1)}+u_{i-1})\\
&=(m+2)\theta_i-\theta_{i-1}.\end{array}\eqno(2.6)$$
Here we are using  the fact that $(a,b)=(a, ax-b)$ for $a, b, x\in \mathbb{Z}$. Thus
$$
\begin{array}{ll}
(u_{p+1}-1,\ u_p)&=(\theta_{-1},\,\theta_0)\\
&=((m+2)\theta_0-\theta_{-1},\theta_0)=(\theta_1,\theta_0)\\
&=(u_p\ ,\ u_{p-1}+1).\end{array}$$
Moreover $(\theta_0, \theta_1)=(\theta_1, (m+2)\theta_1-\theta_0)=(\theta_1, \theta_2)
=\cdots=(\theta_{s-1}, \theta_{s})$, where $s=\lfloor\frac{p}{2}\rfloor$.\\
\indent If $p=2s+1$, then $\theta_{s-1}=u_{s+2}+u_{s-1}=(m+1)h_s$ and $\theta_{s}=u_{s+1}+u_{s}=h_s$.
Thus $$(\theta_{s-1}, \theta_{s})=h_s. $$
\indent If $p=2s$, then
$\theta_{s-1}=u_{s+1}+u_{s-1}=(m+2)u_s$ and
$\theta_s=2u_s$.
Therefore
$$(\theta_{s-1}, \theta_{s})=((m+2)u_s,\,2u_s )=(m,\ 2)u_s. $$
                            $\hfill\Box$\\

\indent The following Lemmas 2.6 and 2.7 will be  used in the proof of Theorem 1.1.\\
\noindent{\bf Lemma 2.6}\, For $n\in\mathbb{N}$, let $B=\begin{pmatrix}
n & \tau_{n-1}   & \tau_{n} \\
0 & \tau_n-\tau_{n-1} & \tau_{n+1}-\tau_n\\
0 & u_n                     & u_{n+1}-1
\end{pmatrix}$
and diag$(s_1(B),\,s_2(B),s_3(B))$  its Smith normal form.\\
\indent If $n=2s+1$, then
$$\left\{\begin{array}{ll}
&s_1(B)=(n,g_s),\\
&s_2(B)=h_s,\\
&s_3(B)=\frac{nh_s}{(n, g_s)}.
\end{array}\right.\eqno(2.7)$$
\indent If $n=2s$, then
$$\left\{\begin{array}{lll}
& s_1(B)=\left(u_s,2\tau_s\right),\\
& s_2(B)=\frac{u_s \left(n,  u_s-4\tau_s\right)}{(u_s,2\tau_s)},\\
& s_3(B)=\frac{n(m+4)u_s}{\left(n,  u_s-4\tau_s\right).}
\end{array}\right.\eqno(2.8)
$$
\noindent{\bf Proof} Recall that $s_1(B)$  equals the greatest common divisor of all entries of $B$. So
 $$
  \begin{array}{ll}
 s_1(B)&=(n,\tau_{n-1}, \tau_n,\tau_n-\tau_{n-1},
 \tau_{n+1}-\tau_n, u_n, u_{n+1}-1)\\
 &=(n,\tau_{n-1}, \tau_n,
 \tau_{n+1}, u_n, u_{n+1}-1). \end{array}$$
Since we have (2.3) and
$$\left\{ \begin{array}{l}
u_n=n-m\tau_n,\\
u_{n+1}-1=(m+1)n+m\tau_{n-1}-m(m+2)\tau_n,
\end{array}\right.$$
it follows that
$$s_1(B)=\left(n, \tau_{n}, \tau_{n-1}\right)
=\left(n, \frac{1}{m}(n-\theta_0),\  \frac{1}{m}(n-\theta_1)\right).\eqno (2.9)$$
From (2.6), we have  $\frac{1}{m}(n-\theta_2)=(m+2) \frac{1}{m}(n-\theta_1)- \frac{1}{m}(n-\theta_0)-n$. Therefore
$$\begin{array}{ll}
s_1(B)&=\left(n, \frac{1}{m}(n-\theta_0),\  \frac{1}{m}(n-\theta_1)\right)\\
&=\left(n, \frac{1}{m}(n-\theta_1),\  \frac{1}{m}(n-\theta_2)\right)\\
&=\cdots \\
&=\left(n, \frac{1}{m}(n-\theta_{s-1}),\  \frac{1}{m}(n-\theta_{s})\right),
  \end{array}
$$ where $s=\lfloor \frac{n}{2} \rfloor$.\\
\indent$\bullet$ If $n=2s+1$, then
$\theta_{s-1}=u_{s+2}+u_{s-1}=(m+1)h_s$ and $\theta_{s}=u_{s+1}+u_{s}=h_s$.
It results that
$$s_1(B)=\left(n,\frac{1}{m}(n-(m+1)h_s), \frac{1}{m}(n-h_s)\right)=\left(n,g_s\right).\eqno(2.10)$$
\indent $\bullet$ If $n=2s$, then
$\theta_{s-1}=u_{s+1}+u_{s-1}=(m+2)u_s$ and $\theta_s=2u_s$.
It results that $$s_1(B)=\left(n,\frac{1}{m}(n-(m+2)u_s),\frac{1}{m}(n-2u_s)
\right)=\left(2s,2\tau_s, u_s\right)=(u_s,2\tau_s).\eqno(2.11)$$
\indent Recall that $s_1(B)s_2(B)$ equals the greatest common divisor of all $2\times 2$ minors of $B$. So
$$s_1(B)s_2(B)=(\Delta_{11}, \Delta_{12}, \Delta_{13}, \Delta_{21}, \Delta_{22}, \Delta_{23},
\Delta_{31}, \Delta_{32}, \Delta_{33}),$$ where $\Delta_{ij}$
is the determinant of the submatrix formed by deleting the $i-$th
row and $j-$th column  of the matrix $B$ for $1\le i,\; j\le 3$.
It is straightforward to see that
$\Delta_{11}=\det\begin{pmatrix}\tau_n-\tau_{n-1} & \tau_{n+1}-\tau_{n} \\
u_n & u_{n+1}-1 \\
\end{pmatrix}
=\frac{1}{m}(u_{n-1}u_{n+1}-u_n^2-1+(u_{n+1}-u_{n-1}))
\overset{(2.4)}{=\!=}\frac{1}{m}(v_n-2);$
$\Delta_{12}=\det\begin{pmatrix}0 &\tau_{n+1}-\tau_n\\
0 &u_{n+1}-1\end{pmatrix}=0$;
$\Delta_{13}=\det\begin{pmatrix}0&\tau_{n}-\tau_{n-1}\\
0& u_n\end{pmatrix}=0;$
$\Delta_{21}
=\det\begin{pmatrix}
          \tau_{n-1} & \tau_n \\
          u_n & u_{n+1}-1
        \end{pmatrix}
=\frac{1}{m}\left((n-1-u_{n-1})(u_{n+1}-1)-(n-u_n)u_{n}\right)
=\frac{1}{m}(n(u_{n+1}-1-u_n))-\frac{1}{m}(u_{n-1}u_{n+1}-u_n^2-1+(u_{n+1}-u_{n-1}))
\overset{(2.4)}{=\!=}n(\tau_n-\tau_{n+1})-\frac{1}{m}(v_n-2);$
$\Delta_{22}=\det\left(\begin{array}{cc}
                     n & \tau_n \\
                     0 & u_{n+1}-1
                   \end{array}
\right)=n(u_{n+1}-1)$;
$\Delta_{23}=\det\begin{pmatrix}
                     n & \tau_{n-1} \\
                     0 & u_{n}
                   \end{pmatrix}=nu_n$;
$\Delta_{31}=\det\begin{pmatrix}
          \tau_{n-1} & \tau_n \\
          \tau_{n}-\tau_{n-1} & \tau_{n+1}-\tau_n
        \end{pmatrix}=\frac{1}{m^2}((n-1-u_{n-1})(u_n-u_{n+1}+1)-(n-u_n)(u_{n-1}+1-u_n))
=\frac{1}{m^2}((u_{n-1}u_{n+1}-u_n^2-1+(-u_{n-1}+u_{n+1}))+2nu_n-n(u_{n+1}+u_{n-1}))
\overset{(2.4)}{=\!=}\frac{1}{m}(\frac{v_n-2}{m}-nu_n);$
$\Delta_{32}=\det\begin{pmatrix}
                     n & \tau_n \\
                     0 & \tau_{n+1}-\tau_n
                   \end{pmatrix}=n(\tau_{n+1}-\tau_n)$;
$\Delta_{33}=\det\begin{pmatrix}
                     n & \tau_{n-1} \\
                     0 & \tau_n-\tau_{n-1}
                   \end{pmatrix}
=n(\tau_n-\tau_{n-1})$. \\
\indent Note that
$\Delta_{33}=\Delta_{23}+\Delta_{32}$,\,
$\Delta_{21}=-\Delta_{32}-\Delta_{11}$ and
$\Delta_{11}=m\Delta_{31}+\Delta_{23}$. So
$$\begin{array}{ll}
s_1(B)s_2(B)&=\left(\Delta_{22},
\Delta_{23},\Delta_{31}, \Delta_{32}\right)\\
&=\left( n(u_{n+1}-1), nu_n,\frac{1}{m}(\frac{v_n-2}{m}-nu_n),n(\tau_{n+1}-\tau_n)\right)\\
&\overset{\theta_i=u_{n-i}+u_i}{=\!=\!=\!=\!=\!=\!=}\left(n\theta_{-1},\, n\theta_0,\,
 \frac{1}{m}(\frac{v_n-2}{m}-n\theta_0),\,\frac{n}{m}(\theta_0-\theta_{-1})\right).\end{array}\eqno(2.12)$$
 \indent With the aid of (2.6), it is easy to verify that
$\frac{n}{m}(\theta_0-\theta_{-1})+n\theta_0=\frac{n}{m}(\theta_1-\theta_0)$.
Moreover, we have
$\frac{1}{m}(\frac{v_n-2}{m}-n\theta_0)-\frac{n}{m}(\theta_1-\theta_0)=\frac{1}{m}(\frac{v_n-2}{m}-n\theta_1)$.
From (2.12), it follows that
$$\begin{array}{ll}
s_1(B)s_2(B)&=\left(n\theta_0, n\theta_1,
\frac{n}{m}(\theta_1-\theta_0),
\frac{1}{m}(\frac{v_n-2}{m}-n\theta_1)\right)\\
&=\left(n\theta_1, n\theta_2, \frac{n}{m}(\theta_2-\theta_1),
\frac{1}{m}(\frac{v_n-2}{m}-n\theta_2)\right)\\
&=\cdots\\
&=\left(n\theta_{s-1}, n\theta_s, \frac{n}{m}(\theta_s-\theta_{s-1}), \frac{1}{m}(\frac{v_n-2}{m}-n\theta_s)\right),
\end{array}\eqno(2.13)$$
where $s=\lfloor\frac{n}{2}\rfloor$.\\
\indent $\bullet$ If $n=2s+1$, then
$\theta_{s-1}=u_{s+2}+u_{s-1}=(m+1)h_s$ and $\theta_s=h_s$. So from (2.4) and (2.13), we can see that
$$\begin{array}{ll}
s_1(B)s_2(B)&=\left(n(m+1)h_s, nh_s, nh_s,  \frac{1}{m}(h_s^2-nh_s)\right)\\
&=\left(nh_s,  \frac{1}{m}(h_s^2-nh_s)\right)\\
&=h_s\left(n, \frac{1}{m}(h_s-n)\right)=h_s(n, g_s).
\end{array}\eqno(2.14)$$
\indent $\bullet$ If  $n=2s$, then
$\theta_{s-1}=u_{s+1}+u_{s-1}=(m+2)u_s$, and  $\theta_s=2u_s$. So from (2.4) and (2.13), we can see that
$$\begin{array}{ll}
s_1(B)s_2(B)&=\left(n\theta_{s-1}, n\theta_s,
\frac{n}{m}(\theta_s-\theta_{s-1}),
\frac{1}{m}((m+4)u_s^2-2nu_s)\right)\\
&=\left(nu_s, \frac{1}{m}((m+4)u_s^2-2nu_s)\right)\\
&=u_s\left(n, \frac{1}{m}((m+4)u_s-2n)\right)\\
&=u_s\left(n, u_s-4\tau_s\right).
\end{array}\eqno(2.15)$$
\indent Recall that $s_1(B)s_2(B)s_3(B)$ equals the determinant of  $B$. So
$$s_1(B)s_2(B)s_3(B)=\det(B)=n\Delta_{11}(B)=\frac{n}{m}(v_n-2).$$
Thus we have
$$\left\{\begin{array}{ll}
s_1(B)s_2(B)s_3(B)=nh_s^2, &\mbox{if}\; n=2s+1,\\
s_1(B)s_2(B)s_3(B)=n(m+4)u_s^2,&\mbox{if}\; n=2s.
\end{array}\right.\eqno(2.16)$$
\indent Combining (2.10), (2.11), (2.14), (2.15), and (2.16),  we obtain the  formulas (2.7) and (2.8).\hfill$\Box$\\

\noindent{\bf Lemma 2.7}\, For $n\in\mathbb{N}$, let $ W=
\begin{pmatrix}
u_{n-1}+1 & u_n\\
u_n     & u_{n+1}-1
\end{pmatrix},$
and  diag$(s_1(W),\,s_2(W))$  its Smith normal form.\\
\indent If $n=2s+1$, then
$$\left\{\begin{array}{ll}
s_1(W)=h_s,\\
s_2(W)=mh_s.\end{array}\right.$$
\indent If $n=2s$, then
$$\left\{\begin{array}{ll}
s_1(W)=(m,2)u_s,\\
s_2(W)=\frac{m(m+4)u_s}{(m, 2)}.\end{array}\right.$$
\noindent{\bf Proof}\, Recall that $s_1(W)$  equals  the greatest common divisor of all the entries of $W$. So
$$\begin{array}{ll}
s_1(W)&=(u_{n-1}+1, u_n, u_{n+1}-1)=(u_{n-1}+1, u_n)\\
&\overset{(2.5)}{=\!=}\left\{\begin{array}{ll}
s_1(W)=h_s, &\mbox{if}\; n=2s+1,\\
s_1(W)=(m,2)u_s, &\mbox{if}\; n=2s.
\end{array}\right.\end{array}\eqno(2.17)$$

\indent  Recall that $s_1(W)s_2(W)$  equals  the greatest common divisor of all $2\times 2$ minors of $W$. So
$$\begin{array}{ll}
s_1(W)s_2(W)&=\det(W)=u_{n-1}u_{n+1}-u_n^2-1+(u_{n+1}-u_{n-1})\\
&\overset{(2.4)}{=\!=}\left\{\begin{array}{ll}
s_1(W)s_2(W)=mh_s^2,& \mbox{if}\; n=2s+1,\\
s_1(W)s_2(W)=m(m+4)u_s^2, & \mbox{if}\; n=2s.
\end{array}\right.\end{array}\eqno(2.18)$$
Combining (2.17) and (2.18), we can obtain
$$s_2(W)=\left\{\begin{array}{cl}
mh_s,                               &\mbox{if}\; n=2s+1,\\
\frac{m(m+4)u_s}{(m, 2)},&\mbox{if}\; n=2s.
\end{array}\right.\eqno(2.19)$$
\hfill$\Box$\\

\section{Proofs of Theorem 1.1 and Corollary 1.2}

\indent  Observe that the critical group of graph $G$ is completely determined by
the cokernel of $L(G)$. Thus, it is sufficient to compute the Smith normal form of the  Laplacian matrix $L(G)$.\\

\indent The proof of Theorem 1.1 contains the following steps:
\begin{enumerate}
\item[(1)] First, we prove that there is a  matrix $A\in \mathbb{Z}^{2m\times 2m}$ such that $L(G) \sim I_{nm-2m} \oplus A$, (see (3.9) and (3.10)).
\item[(2)] Next, we prove there are two matrices $B\in\mathbb{Z}^{3\times 3}$ and $W\in\mathbb{Z}^{2\times 2}$ such that
$A \sim 0_1\oplus B\oplus\underset{m-2}{\underbrace{ W\oplus\cdots\oplus W}}$; the Smith normal forms of $B$ and
$W$ are given in Lemmas 2.6 and 2.7 respectively.
\item[(3)] Finally, we compute the Smith normal form of $A$ from those of $B$ and $W$.
\end{enumerate}
After the three steps, the Smith normal form of $L(G)$ will be obtained.

\subsection*{Step 1}

\indent Now we work on the system of relations of the cokernel of
the Laplacian of $K_m\times C_n$. Let
$e_{i,j}=(0,\cdots,0,1,0,\cdots,0)^t\in \mathbb{Z}^{mn}$, whose unique
nonzero 1 is in the position corresponding to vertex $v_{i,j}$, and
let $x_{i,j}$ be its image in $\coker(L(K_m\times C_n))$. Then it
follows from the relations (1.2) of coker$\left(L(K_m\times C_n)\right)$
that we can get the system of equations:
$$(m+1)x_{i,j}-\sum\limits_{k\in \mathbb{Z}_m\atop k\not=j}x_{i,k}-x_{i-1,j}-x_{i+1,j}=0,\quad
i\in \mathbb{Z}_n,\; j\in \mathbb{Z}_m. \eqno(3.1)
$$
\indent Let $M_i=\sum\limits_{j\in {\mathbb{Z}}_m} x_{i,j}$, for $i\in {\mathbb{Z}}_n$.
Then from (3.1) we have
$$(m+1)M_i-(m-1)M_i-M_{i+1}-M_{i-1}=0.\eqno (3.2)$$
This identity implies that
 $$M_{i+1}=2M_i-M_{i-1}.\eqno(3.3)$$
Recursively using identity (3.3),
we can rewrite all $M_i$'s as integral linear combinations of $M_0$ and $M_1$.
$$M_i=iM_1-(i-1)M_0,\quad 2\le i\le n-1.\eqno(3.4)$$
So from (3.1) and (3.4), we have
$$x_{i,j}=(m+2)x_{i-1,j}-x_{i-2,j}+(i-2)M_0-(i-1)M_1,\quad \eqno(3.5)$$
where $2\le i\le n-1,\quad 0\le j\le m-1.$\\

\noindent{\bf Lemma 3.1}\ For $0\le i\le n-1,\; 0\le j\le m-1$, we have
$$x_{i,j}=-u_{i-1}x_{0,j}+u_ix_{1,j}-\tau_{i-1}M_0+\tau_iM_1. \eqno(3.6)$$
\noindent{\bf Proof}\, This Lemma is valid in cases $i=0,\,1,\,2$.
Suppose that $x_{l,j}=-u_{l-1}x_{0,j}+u_lx_{1,j}-\tau_{l-1}M_0+\tau_lM_1,$ for $l\le
h-1,$ where $h\ge 3$. Then from the induction assumption and the
equations (3.5), it follows that
$$\begin{array}{ll}
x_{h,j}&=(m+2)x_{h-1,j}-x_{h-2,j}+(h-2)M_0-(h-1)M_1,\\
&=(m+2)\left(-u_{h-2}x_{0,j}+u_{h-1}x_{1,j}-\tau_{h-2}M_0+\tau_{h-1}M_1\right)\\
&-\left(-u_{h-3}x_{0,j}+u_{h-2}x_{1,j}-\tau_{h-3}M_0+\tau_{h-2}M_1\right)+(h-2)M_0-(h-1)M_1\\
&=\left(-(m+2)u_{h-2}+u_{h-3}\right)x_{0,j}+((m+2)u_{h-1}-u_{h-2})x_{1,j}\\
&+(-(m+2)\tau_{h-2}+\tau_{h-3}+(h-2))M_0+((m+2)\tau_{h-1}-\tau_{h-2}-(h-1))M_1\\
&=-u_{h-1}x_{0,j}+u_hx_{1,j}-\tau_{h-1}M_0+\tau_hM_1.\end{array}$$
Recall (2.1) and (2.3), we know that (3.6) holds by induction.
\hfill$\Box$\\

\indent In view of Lemma 3.1, we only need  at most
$2m$ generators for the system of equations (3.1). Indeed each $x_{i,j}$ can be expressed in terms of
$x_{0,0},x_{1,0},x_{0,1},x_{1,1},\cdots,\\x_{0,m-1},x_{1,m-1}$. So we know
that there are at least $nm-2m$ diagonal entries of the Smith normal
form  of $L(G)$ are equal to 1 and the remaining invariant  factors
of the abelian group coker$L(K_m\times C_n)$ are the diagonal
entries of the Smith normal form of
the relations matrix induced by  $x_{0,0},x_{1,0},x_{0,1},x_{1,1},\cdots,x_{0,m-1},x_{1,m-1}$.\\

\indent  From (3.6) and the cyclic structure of $K_m\times C_n$, it follows that, for  $0\leq j\leq m-1$,
$$\left\{\begin{array}{ll}x_{0,j}=x_{n,j}=-u_{n-1}x_{0,j}+u_nx_{1,j}-\tau_{n-1}M_0+\tau_nM_1,\\
x_{1,j}=x_{n+1,j}=-u_nx_{0,j}+u_{n+1}x_{1,j}-\tau_nM_0+\tau_{n+1}M_1.\end{array}\right.  \eqno(3.7)$$
Therefore, for $ 0\leq j\leq m-1$,
$$\left\{\begin{array}{ll}
(-u_{n-1}-\tau_{n-1}-1)x_{0,j}+(u_n+\tau_n)x_{1,j}-\tau_{n-1}\sum\limits_{k\not=j}x_{0,k}+\tau_n\sum\limits_{k\not=j}x_{1,k}=0,\\
(-u_n-\tau_n)x_{0,j}+(u_{n+1}+\tau_{n+1}-1)x_{1,j}-\tau_n\sum\limits_{k\not=j}x_{0,k}+\tau_{n+1}\sum\limits_{k\not=j}x_{1,k}=0.
\end{array}\right. \eqno(3.8)$$

\indent Let $$E=\left(
\begin{array}{cc}
-u_{n-1}-1-\tau_{n-1} & u_n+\tau_n\\
-u_n-\tau_n & u_{n+1}-1+\tau_{n+1}
\end{array}
\right),\quad F=\left(\begin{array}{cc}
-\tau_{n-1} & \tau_n\\
-\tau_n & \tau_{n+1}
\end{array}\right),$$
and
$$Y=\left(x_{0,0},x_{1,0},x_{0,1},x_{1,1},\cdots,x_{0,m-1},x_{1,m-1}\right)^T.$$
\indent Then from the equalities in (3.8),  we have that
$$AY=0,\eqno(3.9)$$
where
$$A=\left(
\begin{array}{ccccc}
E & F & F &\cdots & F \\
F & E & F &\cdots & F\\
F & F & E &\cdots & F\\
\vdots &\vdots &\vdots &\ddots &\vdots\\
F & F & F &\cdots & E
\end{array}\right)\in \mathbb{Z}^{2m\times 2m}.\eqno(3.10)$$

\subsection*{Step 2}

\indent The matrix $A$ in equation (3.9) is the relation matrix
induced by the generators
$x_{0,0},x_{1,0},x_{0,1},x_{1,1},\cdots,x_{0,m-1},x_{1,m-1}$. We now discuss
the Smith normal form  of
 the relation matrix $A$.\\
\indent Let $$H=\left(
\begin{array}{ccccc}
I_2       &     0       &    0        & \cdots & 0 \\
-(m-1)I_2 &    I_2      &   I_2       & \cdots & I_2 \\
-I_2      &     0       &   I_2       & \cdots & 0 \\
\vdots    &  \vdots     &  \vdots     & \ddots & 0 \\
-I_2      &     0       &    0        & \cdots & I_2
\end{array}
\right)\in \mathbb{Z}^{2m\times 2m},$$ where $I_2$ is the $2\times 2$ identity matrix. Then it is not difficult to verify that
$$H^{-1}=\left(
\begin{array}{ccccc}
I_2       &      0         &     0         &  \cdots  &  0 \\
I_2       &     I_2        &   -I_2        &  \cdots  & -I_2 \\
I_2       &      0         &    I_2        &  \cdots  &  0 \\
\vdots    &   \vdots       &   \vdots      &  \ddots  &  0 \\
I_2       &      0         &     0         &  \cdots  &  I_2
\end{array}
\right)\in \mathbb{Z}^{2m\times 2m}.$$ By a direct calculation, we have
$$HAH^{-1}= \left(
\begin{array}{ccccc}
    E+(m-1)F     &  F       &      0     &    0     &  0\\
        0        &  E-F     &      0     &    0     &  0\\
        0        &   0      &     E-F    &    0     &  0\\
        0        &   0      &      0     & \ddots   & 0 \\
        0        &   0      &      0     &     0    & E-F
\end{array}\right).\eqno(3.11) $$

\indent Note that $$\begin{pmatrix}E+(m-1)F& F\\
0& E-F\end{pmatrix}= \left(
\begin{array}{cccc}
-n       & n       &   -\tau_{n-1}        & \tau_n       \\
-n       & n       &  -\tau_n                &   \tau_{n+1}     \\
 0       & 0       &      -u_{n-1}-1                      & u_n      \\
 0       & 0       &       -u_n                       &  u_{n+1}-1
\end{array}\right).$$
\indent Let $Q_1= \left(
\begin{array}{cccc}
m  & -m & 1 & -1\\
1  & 0  & 0 & 0 \\
-1 & 1  & 0 & 0 \\
0  & 0  & 0 & 1
\end{array}
\right) $ and\quad  $Q_2=\left(
\begin{array}{cccc}
1 & 0 &   0 & 0 \\
1 & 1 &   0 & 0 \\
0  & 0 &  -1 & 0 \\
0  & 0 &   0 & 1
\end{array}
\right)$.
Then it is easy to see that $Q_1$ and $Q_2$ are unimodular matrices and a careful calculation can show
$$Q_1\begin{pmatrix}E+(m-1)F& F\\
0& E-F\end{pmatrix}Q_2= \left(\begin{array}{cc}
0  & 0 \\
0  & B
\end{array}\right),\eqno(3.12)$$
where the matrix $B$ is just the one defined in Lemma 2.6. \\
\indent Note that $(E-F)
\left(
  \begin{array}{cc}
    -1 & 0 \\
    0 & 1 \\
  \end{array}
\right)
=\begin{pmatrix}
u_{n-1}+1 & u_n\\
u_n     & u_{n+1}-1
\end{pmatrix}=W,$ which is just the one considered in Lemma 2.7. Therefore, from (3.11) and (3.12) we have
$$\begin{array}{ll}
A&\sim HAH^{-1}\sim 0_1\oplus B\oplus\underset{m-2}{\underbrace{ W\oplus\cdots\oplus W}}.\\
&\sim\mbox{diag}(s_1(B), s_2(B), s_3(B), s_1(W), \cdots, s_1(W), s_2(W), \cdots, s_2(W), 0)\\
&\sim\mbox{diag}(s_1(B), s_2(B), s_1(W), \cdots, s_1(W), s_2(W),
\cdots, s_2(W), s_3(B), 0).
\end{array}\eqno(3.13)$$

\subsection*{Step 3}

\indent Now we distinguish two cases to compute the Smith normal form of $A$.\\

\indent $Case\, 1.$\, $n=2s+1$.\\
\indent Then by Lemmas 2.6 and 2.7 we have that   $s_1(B)=(n, g_s), s_2(B)=h_s,
 s_3(B)=\frac{nh_s}{(n, g_s)},
 s_1(W)=h_s, s_2(W)=mh_s.$
Since $s_1(B)|s_2(B)$, $s_1(W)|s_2(W)$ and $s_2(B)=s_1(W)$, it follows that $s_1(B)|s_2(W)$.\\
\indent Write $\gamma=(s_2(W), s_3(B))$ and $\varphi=\frac{s_2(W)s_3(B)}{\gamma}$. Then
$$\mbox{diag}(s_2(W), s_3(B))\sim\mbox{diag}(\gamma, \varphi),\eqno(3.14) $$
where
$$\begin{array}{ll}
\gamma &=\left(mh_s, \frac{nh_s}{(n, g_s)}\right)=\frac{h_s}{(n, g_s)}(m(n, g_s), n)\\
&=\frac{h_s}{(n, g_s)}(mn, n-h_s, n)\quad(\mbox{Here}\;  mg_s=n-h_s.)\\
&=\frac{h_s}{(n, g_s)}(n, h_s),
\end{array}\eqno(3.15)$$
and $$\varphi=\frac{nmh_s}{(n, h_s)}.\eqno(3.16)$$

Note that $s_1(W)=s_2(B)$, it implies that $s_1(W)|s_3(B)$. So
$s_1(W)|\gamma|s_2(W)$. Moreover, it is easy to see that $s_2(W)\mid\varphi$.\\
\indent Therefore, (3.13) and (1.14) implies that
$$\mbox{diag}(s_1(B), s_2(B), \underset{m-2}{\underbrace{s_1(W), \cdots,
s_1(W)}}, \gamma, \underset{m-3}{\underbrace{s_2(W), \cdots,
s_2(W)}}, \varphi, 0)\eqno(3.17)$$
is the Smith normal form of $A$.\\

\indent $Case\, 2.$\, $n=2s$.\\
\indent Then by Lemmas 2.6 and 2.7, we know that in this case $s_1(B)=(u_s, 2\tau_s)$, $s_2(B)=\frac{u_s
(n,  u_s-4\tau_s)}{(u_s, 2\tau_s)}$,
$s_3(B)=\frac{n(m+4)u_s}{(n,  u_s-4\tau_s)}$,
$s_1(W)=(m, 2)u_s$,
$s_2(W)=\frac{m(m+4)u_s}{(m, 2)}.$\\
\indent It is obvious that $$s_1(B)|s_2(B)|s_3(B),\ s_1(B)|s_1(W)|s_2(W), \ s_2(B)|s_2(W). \eqno(3.18)$$
\indent Clearly, we have
$$\mbox{diag}(s_2(W), s_3(B))\sim\mbox{diag}(\rho, \xi), \eqno(3.19)$$
where $$\begin{array}{lll}\rho &=(s_2(W), s_3(B))=\left(\frac{m(m+4)u_s}{(m, 2)},\frac{n(m+4)u_s}
{(n,  u_s-4\tau_s)}\right)\\
&=\frac{(m+4)u_s}{(n,\, u_s-4\tau_s)(m,\, 2)}( mn,\, (m+4)u_s-2n,\, mn,\, 2n)\\
&=\frac{(m+4)u_s( mn,\, (m+4)u_s,\, 2n)}{(n,\, u_s-4\tau_s)(m,\, 2)},\end{array}\eqno(3.20)$$
and
$$\xi=\frac{s_2(W)s_3(B)}{\rho}=\frac{mn(m+4)u_s}{\left(mn,\ (m+4)u_s,\ 2n\right)}.\eqno(3.21)$$
Therefore,  It follows from (3.13) that
$$A\sim \mbox{diag}(s_1(B), s_2(B), \overset{m-2}{\overbrace{s_1(W), \cdots, s_1(W)}}, \rho,
\overset{m-3}{\overbrace{s_2(W), \cdots, s_2(W)}}, \xi, 0).\eqno(3.22)$$
\indent We also have
$$\mbox{diag}(s_2(B), s_1(W))\sim\mbox{diag}(\zeta, \eta),\eqno(3.23)$$
where
$$\begin{array}{lllll}
\zeta
&=(s_2(B), s_1(W))=\left( \frac{u_s(n,  u_s-4\tau_s)}{(u_s, 2\tau_s)},(m, 2)u_s\right)\\
&=\frac{u_s\left(n,\, u_s-4\tau_s,\, mu_s,\, 2m\tau_s,\, 2u_s,\, 4\tau_s\right)}{(u_s, 2\tau_s)}\\
&=\frac{u_s\left(n,\, u_s,\, 4\tau_s\right)}{(u_s, 2\tau_s)},
\end{array}\eqno(3.24)$$
and
$$\eta=\frac{s_2(B)s_1(W)}{\zeta}=\frac{u_s(m,2) \left(n,  u_s-4\tau_s\right)}
{\left(n,\, u_s,\, 4\tau_s\right)}.\eqno(3.25)$$

Combining (3.22) and (3.23), we have
$$A\sim \mbox{diag}(s_1(B), \zeta,\underset{m-3}{\underbrace{s_1(W),\cdots,
s_1(W)}}, \eta, \rho,\underset{m-3}{\underbrace{s_2(W),\cdots,
s_2(W)}}, \xi, 0).$$

\indent Since $\left(n,  u_s-4\tau_s\right)|n$, and $\left(m,
2\right)u_s|(m+4)u_s$, it follows that $s_1(W)|s_3(B)$ and hence
$s_1(W)|(s_2(W), s_3(B))$. Furthermore, since $s_2(B)|s_3(B)$ and
$s_2(B)|s_2(W)$, it implies that $s_2(B)|(s_2(W), s_3(B))$.
 So, $(s_2(W), s_3(B))$, i.e., $\rho$
is a common multiple of $s_1(W)$ and $s_2(B)$. Note that $\eta$ is
the least common multiple of $s_1(W)$ and $s_2(B)$, it divides
$\rho$.  According to $(3.18)$, it is easy to see $s_1(B)|\zeta$. And
it is clear that we have $s_1(W)|\eta$ and $s_2(W)|\xi$. Thus
$$\mbox{diag}(s_1(B), \zeta,\underset{m-3}{\underbrace{s_1(W),\cdots,
s_1(W)}}, \eta, \rho,\underset{m-3}{\underbrace{s_2(W),\cdots,
s_2(W)}}, \xi, 0)\eqno(3.26)$$
is the Smith normal form of $A$.\\
\indent Now, the proof of Theorem 1.1 is completed. \hfill$\Box$\\

\noindent{\bf Proof of Corollary 1.2}\, If $n=2s+1$, then
$s_1(W)s_2(W)=mh_s^2\overset{(2.4)}{=\!=\!=}v_n-2$ and
$s_1(B)s_2(B)s_1(W)\gamma\varphi=(n,g_s)h_sh_s\frac{h_s}{(n,g_s)}(n,
h_s)
\frac{nmh_s}{(n,h_s)}=nmh_s^4\overset{(2.4)}{=\!=\!=}\frac{n}{m}(v_n-2)^2$.
It follows that  the spanning tree number of $K_m\times C_n$ is
$\frac{n}{m}(v_n-2)^2\times(v_n-2)^{m-3}=\frac{n}{m}(v_n-2)^{m-1}$.\\
\indent If $n=2s$, then $s_1(W)s_2(W)=\left(m,
2\right)u_s\frac{m(m+4)u_s}{\left(m,
2\right)}=m(m+4)u_s^2\overset{(2.4)}{=\!=\!=}v_n-2$ and
$s_1(B)\zeta\eta\rho\xi=nm(m+4)^2u_s^4\overset{(2.4)}{=\!=\!=}\frac{n}{m}(v_n-2)^2$. It
follows that  the spanning tree number of $K_m\times C_n$ is
$\frac{n}{m}(v_n-2)^2\times(v_n-2)^{m-3}=\frac{n}{m}(v_n-2)^{m-1}$.\\

\noindent (In fact, by (3.13) we know that the spanning tree number of graph
$K_m\times C_n$ is $|\det( \mbox{diag}(B\oplus \underbrace{W\oplus\cdots\oplus W}_{m-2}))|=\det (B)\times (\det (W))^{m-2}$.
In the proofs of Lemmas 2.6 and 2.7, we have seen that $\det (B)=\frac{n}{m}(v_n-2)$ and $\det (W)=v_n-2$.  So the spanning tree number
is $\frac{n}{m}(v_n-2)^{m-1}$.)\hfill$\Box$\\

\section{Remarks}

\indent\, (\uppercase\expandafter {\romannumeral 1})\; If $n=1$,  then $K_m\times C_n$ is the complete
graph $K_m$, from [12] we know its critical group is $(\mathbb{Z}_m)^{m-2}$.\\
\indent\, (\uppercase\expandafter {\romannumeral 2})\;
 If $n=2$, then $\coker(K_m\times C_2)$ is determined by the generators
 $x_{0,j}$, $x_{1,j}$, and the relations
$$\hspace{3.2cm}\left\{\begin{array}{ll}
(m+1)x_{0,j}-\sum\limits_{k\in\mathbb{Z}_m\atop k\not=j} x_{0,k}-2x_{1,j}=0,\hspace{4cm} (4.1)\\
\\
(m+1)x_{1,j}-\sum\limits_{k\in\mathbb{Z}_m\atop k\not=j} x_{1,k}-2x_{0,j}=0,\hspace{4cm} (4.2)
\end{array}\right.$$
Where $j\in\mathbb{Z}_m$.\; From (4.1), we get
$$2x_{1,j}=(m+1)x_{0,j}-\sum\limits_{k\in\mathbb{Z}_m\atop k\not=j}  x_{0,k},\quad  j\in\mathbb{Z}_m.\eqno(4.3)$$
Substituting (4.3) into $(4.2)\times 2$ gives the following
$$(m+4)(m-1)x_{0,j}-\sum\limits_{k\in\mathbb{Z}_m\atop k\neq j} (m+4)x_{0,k}=0,\quad j\in\mathbb{Z}_m.$$
\indent So we can simplify  (up to equivalence) the Laplacian matrix of $K_m\times C_2$ into
$$(m+4)\left(
    \begin{array}{cccc}
      m-1 & -1 & \cdots & -1  \\
      -1 & m-1 & \cdots & \vdots  \\
      \vdots & \vdots & \ddots & -1  \\
      -1 & \cdots & -1 & m-1  \\
        \end{array}
  \right)_{m\times m}=(m+4)L(K_m).$$
Therefore,  the critical group of $K_m\times C_2$ is
$\mathbb{Z}_{m+4}\oplus(\mathbb{Z}_{m(m+4)})^{m-2}$.\\

\indent (\uppercase\expandafter {\romannumeral 3})
 If $m=1$, then  $K_m\times C_n$ is just the
cycle $C_n$. So from [14], its critical group is $\mathbb{Z}_n$. In fact,
when $m=1$, it is easy to see the matrix $A$ in (3.10) has the
following property:
$$A=E=\left(\begin{array}{cc}
-n & n\\
-n & n
\end{array}\right)
\sim (0)\oplus n.$$
The known result is obtained immediately.\\
\indent (\uppercase\expandafter {\romannumeral 4}) If $m=2$, the graph $K_2\times C_n$ is just the
Cayley graph $\mathcal{D}_n$ of dihedral group. The result of this
case was obtained  in [8].
In the following  we will try to get the result again.\\
\indent From (3.10) and (3.13), we know $$A=\left(\begin{array}{cc}
E & F\\
F & E
\end{array}\right)\sim(0)\oplus B
\sim\mbox{diag}(s_1(B), s_2(B), s_3(B), 0).$$

\indent $\bullet$ If $n=2s+1$, then from (2.7) we have
$$\left\{\begin{array}{ll}
&s_1(B)=(n,g_s),\\
&s_2(B)=h_s,\\
&s_3(B)=\frac{nh_s}{(n, g_s)}.
\end{array}\right.$$
\indent Note that $g_s=\frac{n-h_s}{2}$ and $n$ is odd, then
$(n,g_s)=(n,2g_s)=(n,n-h_s)=(n,h_s)$. \indent Therefore,
$$\left\{\begin{array}{ll}
&s_1(B)=(n,h_s),\\
&s_2(B)=h_s,\\
&s_3(B)=\frac{nh_s}{(n, h_s)}.
\end{array}\right.$$

\indent $\bullet$ If $n=2s$, then from (2.8) we have
$$\left\{\begin{array}{lll}
& s_1(B)=(u_s, 2\tau_s),\\
& s_2(B)=\frac{u_s(n, u_s-4\tau_s)}{(u_s, 2\tau_s)},\\
& s_3(B)=\frac{6nu_s}{(n,  u_s-4\tau_s)}.
\end{array}\right.
$$
Note that $s-u_s$ is even(2.2), and $2^{t+1}|u_n$ if $2^t|n$
(Corollary 3.4 [8]). Hence
$$\begin{array}{ll}
s_1(B)&=(u_s,2\tau_s)=(u_s,s-u_s)=(u_s, s)\\
&=\left\{\begin{array}{ll}
(u_s,2s)=(u_s,n), &  \mbox{if $s$ is odd,}\\
\frac{(u_s,2s)}{2}=\frac{(u_s,n)}{2}, &
\mbox{if $s$ is even.}
\end{array}\right.
\end{array}$$
Also we have $3^t|u_n$ if $3^t|n$ (Corollary 3.4 [8]), then
$$(n, u_s-2(s-u_s))=(2s,3u_s)=(n,u_s).$$
So $$s_2(B)=\frac{u_s(n, u_s-4\tau_s)}{(u_s,2\tau_s)}=
\frac{u_s(n,u_s)}{(u_s,2\tau_s)}=\left\{
\begin{array}{ll}
u_s, &  \mbox{if  $s$ is  odd,}\\
2u_s, & \mbox{if $s$ is even.}
\end{array}\right.$$
Then $s_3(B)=\frac{6nu_s}{(n,  u_s-4\tau_s)}=\frac{6nu_s}{(n,u_s)}$.\\

\indent It is easy to see the result of case $m=2$ here is the same
to the result obtained in [8].

\section*{Acknowledgements}

\indent Many thanks to the  referee for his/her many helpful comments and
 suggestions, which have considerably improved the presentation of this paper.

\end{document}